\documentclass[11pt]{amsart}\usepackage{}\usepackage{cases}\usepackage{amsmath}\usepackage{amssymb}\usepackage{amsfonts}
\newcommand{\Pf}{\noindent {\bf Proof }}
\newcommand{\Rk}{\noindent {\bf Remark}}
\newtheorem{thm}{Theorem}[section]
\newtheorem{prop}{Proposition}[section]
\newtheorem{lemma}{Lemma}[section]
\newtheorem{coro}{Corollary}[section]

\begin{document}
\title[Weitzenb\"ock formula for canonical metrics]{Berger Curvature Decomposition, Weitzenb\"ock Formula, and Canonical Metrics on Four-manifolds}
\author{Peng Wu}
\address{Department of Mathematics, Cornell University, Ithaca, NY 14853, United States} \email{wupenguin@math.cornell.edu}
\thanks{}
\dedicatory{}
\date{May 20, 2014}


\begin{abstract}
We first provide an alternative proof of the classical Weitzneb\"ock formula for Einstein four-manifolds using Berger curvature decomposition, motivated by which we established a unified framework for a Weitzenb\"ock formula for a large class of canonical metrics on four-manifolds (or a Weitzenb\"ock formula for ``Einstein metrics" on four-dimensional smooth metric measure spaces). As applications, we classify Einstein four-manifolds of half two-nonnegative curvature operator, which in some sense provides a characterization of K\"ahler-Einstein metrics with positive scalar curvature on four-manifolds, we also discuss four-manifolds of half two-nonnegative curvature operator and half harmonic Weyl curvature.
\end{abstract}
\maketitle

\section{Introduction}

For an oriented Riemannian four-manifold $(M,\ g)$, the Hodge star operator $\star: \wedge^2 TM\rightarrow \wedge^2 TM$ induces an eigenspace decomposition $\wedge^2 TM=\wedge^+ M\oplus \wedge^- M$, where $\wedge^{\pm} M=\{\omega\in\wedge^2TM: \star \omega=\pm\omega\}$ are eigenspaces of the Hodge star operator. Elements in $\wedge^{\pm} M$ are called self-dual and anti-self-dual 2-forms. This decomposition further induces the duality decomposition of the curvature operator $\mathfrak{R}:\wedge^2 TM\rightarrow \wedge^2 TM$,
\begin{equation}
\label{dual} \mathfrak{R} =\left( \begin{array}{cc}
\frac{R}{12}g+W^+ & \overset{\circ}{\mathrm{Ric}} \\
\overset{\circ}{\mathrm{Ric}} & \frac{R}{12}g+W^-
\end{array} \right),
\end{equation}
where $R$ is the scalar curvature, $W^{\pm}$ are called self-dual and anti-self-dual parts of Weyl curvature tensor. $(M,\ g)$ is called (anti-)self-dual, or half conformally flat, if $W^-=0$ ($W^+=0$).

For Einstein four-manifolds, $\mathfrak{R}$, $W$, $\mathfrak{R}^{\pm}$, and $W^{\pm}$ are all harmonic. Using the harmonicity, Derdzi\'nski \cite{Der} derived the following Weitzenb\"ock formula,
\begin{thm} [\cite{Der}] \label{FourWeitzenbock}
Let $(M, g)$ be an oriented Einstein four-manifold with $\mathrm{Ric}=\lambda g$, then
\begin{equation*}
\begin{split}
\Delta |W^{\pm}|^2 =& 2|\nabla W^{\pm}|^2+4\lambda |W^{\pm}|^2-36\det W^{\pm}.
\end{split}
\end{equation*}
where $\langle S, T\rangle=\frac{1}{4}S_{ijkl}T^{ijkl}$ for any $(0,4)$-tensor $S, T$.
\end{thm}

The Weitzenb\"ock formula, together with Hitchin's classification of half conformally flat Einstein four-manifolds \cite{Besse}, plays a key role in the classification of Einstein four-manifolds of positive scalar curvature, see for example \cite{GL, Yang}. In \cite{GL}, Gursky, LeBrun obtained an optimal gap theorem for $W^{\pm}$ and classified Einstein four-manifolds with nonnegative sectional curvature operator and positive intersection form. In \cite{Yang}, Yang classified Einstein four-manifolds with $\text{Ric}=g$ and sectional curvature bounded below by $\frac{\sqrt{1249}-23}{120}$. 
\vspace{0.4cm}

In this paper, first following from an argument in the author's Ph.D. thesis \cite{Wuthesis}, we provide an alternative proof of the Weitzenb\"ock formula in Theorem \ref{FourWeitzenbock} by combining an argument of Hamilton (Lemma 7.2 in \cite{Hamilton}) and Berger curvature decomposition \cite{Berger}. As an application, we classify Einstein four-manifolds of half two-nonnegative curvature operator (half nonnegative isotropic curvature).

A Riemannian metric is said to have $k$-positive ($k$-nonnegative) curvature operator if the sum of any $k$ eigenvalues is positive (nonnegative). A Riemannian metric on a four-manifold is said to have half two-positive (two-nonnegative) curvature operator, 
if the self-dual curvature operator $\mathfrak{R}^+=\frac{R}{12}g+W^+$ or the anti-self-dual curvature operator $\mathfrak{R}^-=\frac{R}{12}g+W^-$ is two-positive (two-nonnegative). By the duality decomposition, it is easy to check that half two-positive curvature operator and half positive isotropic curvature are equivalent. It is obvious that if a four-manifold $(M,\ g)$ is half conformally flat and has positive scalar curvature, then 
$\mathfrak{R}$ is half two-positive. Another interesting fact is that any K\"ahler metric on a four-manifold has two-nonnegative $\mathfrak{R}^+$.

Precisely, we prove

\begin{thm} \label{half}
Let $(M,\ g)$ be a simply-connected Einstein four-manifold with positive scalar curvature.

(1). If $\mathfrak{R}^+$ or $\mathfrak{R}^-$ is two-positive, then it is half conformally flat, hence isometric to $(S^4, g_0)$ or $(\mathbb{C}P^2, g_{FS})$.

(2). If $\mathfrak{R}^+$ or $\mathfrak{R}^-$ is two-nonnegative, then it is isometric to $(S^4, g_0)$ or a K\"ahler-Einstein surface.

(3). If $\mathfrak{R}^+$ or $\mathfrak{R}^-$ is two-nonnegative, and $\mathcal{R}$ is four-nonnegative, then it is isometric to $(S^4, g_0)$, $(\mathbb{C}P^2, g_{FS})$, or $(S^2\times S^2,\ g_0\oplus g_0)$.

\end{thm}

\Rk: After the author finished the proof of Theorem \ref{half}, he observed that it was proved independently by Richard and Seshadri \cite{RS} using a different method. They first proved that the cone of half nonnegative isotropic curvature is preserved along the Ricci flow, then applied an argument of Brendle in \cite{Brendle}.

Recall that for any K\"ahler metric on a four-manifold, $W^+$ has eigenvalues $\{\frac{R}{6},-\frac{R}{12},-\frac{R}{12}\}$, hence any K\"ahler metric has half two-nonnegative curvature operator (or half nonnegative isotropic curvature). Derdzinski proved that if a Riemannian metric on a four-manifold satisfies $\delta W^+=0$ and $W^+$ has at most two distinct eigenvalues, then the metric is locally conformally K\"ahler, if in addition the scalar curvature is constant, then the metric itself is K\"ahler.
It is interesting to point out that part (2) of Theorem \ref{half} in fact provides a characterization of K\"ahler-Einstein metrics of positive scalar curvature on four-manifolds: any Einstein metric which is not conformally flat and has half two-nonnegative curvature operator is K\"ahler-Einstein.

\

Motivated by our alternative proof of Weitzenb\"ock formula for Einstein four-manifolds, we establish a unified framework for the Weitzenb\"ock formula for a large class of metrics on four-manifolds, which are called generalized $m$-quasi-Einstein metrics.

Let $(M^n, g)$ be a Riemannian manifold, $g$ is called a \textit{generalized m-quasi-Einstein metric} \cite{Catino2}, if
\begin{equation} \label{quasiEinstein}
\text{Ric}_f^m=\text{Ric}+\nabla^2 f-\frac{1}{m}df\otimes df=\lambda g,
\end{equation}
for some $f,\lambda\in C^{\infty}(M)$ and $m\in\mathbb{R}\cup\{\pm\infty\}$. Notice that $\mathrm{Ric}_f^m$ is exactly the $m$-Bakry-Emery Ricci curvature, introduced by Bakry and Emery \cite{BE}, for smooth metric measure spaces, therefore a generalized $m$-quasi-Einstein metric on a Riemannian manifold can be considered as an ``Einstein metric" on a smooth metric measure space. In particular, ``Einstein metrics" on smooth metric measure spaces contain at least the following interesting special cases,

(1), when $f=$const, it is an Einstein metric;

(2), when $m=\infty$ and $\lambda=$const, it is a gradient Ricci soliton;

(3), when $0<m<\infty$ and $\lambda=$const, it is an $m$-quasi-Einstein metric, and $(M^n\times F^m, g_M+e^{-\frac{2f}{m}}g_F)$ is a warped product Einstein manifold, where $F^m$ is an $m$-dimensional space form;

(4), when $m=1$, together with other conditions, it is a static metric in general relativity;

(5), when $m=2-n$, it is a conformally Einstein metric, and $\bar g=e^{\frac{2}{2-n}f}g$ is an Einstein metric. 
\vspace{0.4cm}

The Weitzenb\"ock formual can be stated as following,

\begin{thm} \label{UnifiedWeitzenbock}
Let $(M, g)$ be a generalized $m$-quasi-Einstein four-manifold with $\mathrm{Ric}_f^m=\lambda g$, then
\begin{equation*}
\begin{split}
\Delta_f |W^{\pm}|^2=2|\nabla W^{\pm}|^2+&(4\lambda+\frac{2}{m}|\nabla f|^2)|W^{\pm}|^2-36\det W^{\pm}\\
-&(1+\frac{2}{m})\langle(\nabla^2f\circ\nabla^2f)^{\pm},W^{\pm}\rangle,
\end{split}
\end{equation*}
where $T^{\pm}(\alpha,\beta)=T(\alpha^{\pm}, \beta^{\pm})$ for $\alpha,\beta\in\wedge^2 M$.
\end{thm}

As special cases, we get the Weitzenb\"ock formula for conformally Einstein four-manifolds and four-dimensional gradient Ricci
solitons,

\begin{coro} \label{WconfEinstein}
Let $(M,g,f)$ be a conformally Einstein four-manifolds with $\mathrm{Ric}+\nabla^2f+\frac{1}{2}df\otimes df=\lambda g$. Then
\begin{equation*}
\begin{split}
\Delta_f |W^{\pm}|^2=& 2|\nabla W^{\pm}|^2+(4\lambda-|\nabla f|^2)|W^{\pm}|^2-36\det W^{\pm}.
\end{split}
\end{equation*}
\end{coro}

\begin{coro} \label{WSoliton}
Let $(M, g, f)$ be a four-dimensional gradient Ricci soliton with $\mathrm{Ric}+\nabla^2f=\lambda g$. Then
\begin{equation*}
\begin{split}
\Delta_f |W^{\pm}|^2=& 2|\nabla W^{\pm}|^2+4\lambda|W^{\pm}|^2-36\det W^{\pm} -\langle(\mathrm{Ric}\circ\mathrm{Ric})^{\pm},W^{\pm}\rangle,
\end{split}
\end{equation*}
\end{coro}

There have been several generalizations of the Weitzenb\"ock formula for Einstein four-manifolds.

(1). In \cite{CGY}, Chang, Gursky, and Yang derived an integral Weitzenb\"ock formula for all compact four-manifolds,
\begin{equation*} \label{intWeitzenbock}
\begin{split}
0=\int_M2|\nabla W^{\pm}|^2-8|\delta W^{\pm}|^2+R|W^{\pm}|^2-36\det W^{\pm},
\end{split}
\end{equation*}
and they also derived an integral Weitzenb\"ock formula for Bach-flat metrics, with the help of which they proved a very interesting conformally invariant sphere theorem in four dimensions.

(2). The Weitzenb\"ock formula for conformally Einstein four-manifolds in Corollary \ref{WconfEinstein} can also be derived directly using the property of the conformal change of $\delta W^{\pm}$ and the Weitzenb\"ock formula for Einstein four-manifolds, see for example \cite{Der,Gursky,LeBrun}. Since gradient Ricci solitons are self-similar solutions to the Ricci flow, the Weitzenb\"ock formula for gradient Ricci solitons in Corollary \ref{WSoliton} can also be derived by applying the observation of the author \cite{Wuthesis} and the evolution equation of the Weyl curvature in \cite{CaMa}, see \cite{CT}.

\

We observe from Theorem \ref{FourWeitzenbock} that,

\begin{thm} \label{harmonicWeyl}
Let $(M,g)$ be a compact four-dimensional Riemannian manifold. If $\delta W^{\pm}=0$ and $\mathfrak{R}^{\pm}$ is two-positive, then $g$ is either self-dual or anti-self-dual. If $\delta W^{\pm}=0$ and $\mathfrak{R}^{\pm}$ is two-nonnegative, then either $g$ is self-dual or anti-self-dual, or $g$ is a cscK metric.
\end{thm}

The proof is based on an observation for half two-nonnegative curvature operator, see Lemma \ref{halfPIC} in Section 3. If in addition that $g$ is a gradient Ricci soliton on $M$, then we get

\begin{coro}
Let $(M,g,f)$ be a compact four-dimensional gradient shrinking Ricci soliton. If $\delta W^{\pm}=0$ and $\mathfrak{R}^{\pm}$ is two-nonnegative, then $(M,g)$ is isometric to $(S^4,g_0)$ or K\"ahler-Einstein.
\end{coro}

Gradient Ricci solitons were introduced by Hamilton \cite{Hamilton2}, they played an important role in the Ricci flow and Perelman's resolution to the Poincar\'e conjecture and the geometrization conjecture \cite{Perel1, Perel2, Perel3}. In the past three decades, there has been lots of work on the classification of gradient shrinking Ricci solitons. In dimensions 2 and 3, by \cite{Hamilton2, Ivey, Perel1, NW, CCZ}, the classification is complete. In dimensions greater than or equal to 4, by \cite{ENM,NW,Zhang,PW,CWZ,MS,FG,CaoChen} and references therein, the classifications of gradient shrinking Ricci solitons with vanishing Bach tensor or harmonic Weyl curvature are complete. In particular in dimension 4, half conformally flat gradient shrinking Ricci
solitons have been completely classified in \cite{CaoChen,ChenWang}, which can be considered as an analogue of Hitchin's calssical classification of half conformally flat Einstein four-manifolds.

Further applications of the Weitzenbock formula to gradient Ricci solitons and conformally Einstein four-manifolds will be addressed in subsequent work \cite{CW, WWW}.
\vspace{0.3cm}

The paper is organized as following. In Section 2 we discuss Berger curvature decomposition, provide an alternative proof of the Weitzenb\"ock formula for Einstein four-manifolds using Berger curvature decomposition. and classify Einstein four-manifolds of half two-nonnegative curvature operator. In Section 3, we prove the Weitzenb\"ock formula for generalized quasi-Einstein manifolds and Theorem \ref{harmonicWeyl}. In the appendix we provide Berger's proof of Berger curvature decomposition.

\textbf{Acknowledgement}. This paper is another extension of the author's Ph.D. thesis, he expresses his great attitude to his advisors Professors Xianzhe Dai and Guofang Wei for their guidance, encouragement, and constant support. He thanks Professor Zhenlei Zhang for bringing the classification of Einstein four-manifolds of half two-positive curvature operator to the author's attention. He thanks Professor Jeffrey Case for helpful discussions.
The first part of the work was done when the author was visiting BICMR in summer 2013. He thanks the institute for their hospitality and Professor Yuguang Shi for his help. The author was partially supported by an AMS-Simons postdoctoral travel grant.


\section{Berger Curvature Decomposition and Proof of Theorem 1.1 and 1.2 }

First we fix the notations. Our sign conventions for the curvature tensor will be so that

\begin{equation*}
\begin{split}
&R_{ijkl}=g_{hk}R^h_{ijl},\ \ K(e_i,e_j)=R_{ijij},\ \
R_{ik}=g^{jl}R_{ijkl},\ \ R=g^{ij}R_{ij}.\\
&(\nabla_p\nabla_q-\nabla_q\nabla_p)T_{ijkl}=R_{pqim}T_{mjkl}+...+R_{pqlm}T_{ijkm}.\\
\end{split}
\end{equation*}
And our convention for the inner product of two $(0,4)$-tensors $S,T$ will be
\begin{equation*}
\langle S,T\rangle=\frac{1}{4}S_{ijkl}T^{ijkl}\
\end{equation*}
so that our convention agrees with the  one in Derdzinski's Weitzenb\"ock formula \cite{Der}.

We start from an interesting observation of Berger \cite{Berger},

\begin{lemma} \label{BergerLemma}
Let $(M, g)$ be an oriented Einstein four-manifold. Then for any $p\in M$, and any orthonormal basis $\{e_1, e_2, e_3, e_4\}$ of $T_p M$,
\begin{equation*}
\begin{split}
K(e_1, e_2)&= K(e_3,e_4),\\
K(e_1, e_3)&= K(e_2,e_4),\\
K(e_1, e_4)&= K(e_2,e_3).
\end{split}
\end{equation*}
\end{lemma}

In another word, Lemma \ref{BergerLemma} says that for an Einstein four-manifold,
\begin{equation*}
\begin{split}
R_{ijkl}=&R_{i'j'k'l'},
\end{split}
\end{equation*}
where $(i'j')$ is the dual of the pair $(ij)$, i.e., the pair such that $e_i\wedge e_j\pm e_{i'}\wedge e_{j'}\in \wedge^{\pm}M$, or $(iji'j')=\sigma(1234)$ for some even permutation $\sigma\in S_4$.

Using Lemma \ref{BergerLemma} and basic symmetries of curvature tensor, Berger obtained the following curvature decomposition \cite{Berger} for Einstein four-manifolds (see also \cite{ST}), see the appendix for the proof,

\begin{prop} \label{Bergerdecomposition}
Let $(M,\ g)$ be an Einstein four-manifold with
$\mathrm{Ric}=\lambda g$. Then at any $p\in M$, there exists an orthonormal basis $\{e_i\}_{1\leq i\leq 4}$ of $T_p M$, such that relative to the corresponding basis $\{e_i\wedge e_j\}_{1\leq i<j\leq 4}$ of $\wedge^2 T_pM$, $\mathfrak{R}$ takes the form
\begin{equation*}
\mathfrak{R}=\left( \begin{array}{cc}
A & B\\
B & A
\end{array}\right),
\end{equation*}
where $A=\text{diag}\{a_1,\ a_2,\ a_3\}$, $B=\text{diag}\{b_1,\ b_2,\ b_3\}$ satisfy the following properties,\\
$\mathrm{(1)}$. $a_1=K(e_1, e_2)=K(e_3, e_4)=\min\{K(\sigma): \sigma\in\wedge^2 T_p M,||\sigma||=1\}$,

\hspace{0.22cm} $a_3=K(e_1,e_4)=K(e_2,e_3)=\max\{K(\sigma):\sigma\in\wedge^2 T_p M, ||\sigma||=1\}$,

\hspace{0.22cm} $a_2=K(e_1, e_3)=K(e_2, e_4)$, $a_1+a_2+a_3=\lambda$;\\
$\mathrm{(2)}$. $b_1=R_{1234},\ b_2=R_{1342},\ b_3=R_{1423}$;\\
$\mathrm{(3)}$. $|b_2-b_1|\leq a_2-a_1,\ |b_3-b_1|\leq a_3-a_1,\ |b_3-b_2|\leq a_3-a_2$.
\end{prop}

As observed in \cite{Wuthesis}, Berger curvature decomposition is in fact a special case of the duality decomposition, as it is easy to see that the eigenvalues of $\mathfrak{R}$ are
\begin{equation*}
\begin{split}
& a_1+b_1\leq a_2+b_2\leq a_3+b_3,\\
& a_1-b_1\leq a_2-b_2\leq a_3-b_3,
\end{split}
\end{equation*}
with corresponding eigenvectors $\omega_1^{\pm}=\frac{1}{\sqrt{2}}(e_1\wedge e_2 \pm e_3\wedge e_4),\ \omega_2^{\pm}=\frac{1}{\sqrt{2}}(e_1\wedge e_3 \pm e_4\wedge e_2),\ \omega_3^{\pm}=\frac{1}{\sqrt{2}}(e_1\wedge e_4 \pm e_2\wedge e_3)$. In another word, for Einstein four-manifolds,
\begin{equation}
\label{dual} \mathfrak{R} =\left( \begin{array}{cc}
\frac{R}{12}g+W^+ & 0\\
0 & \frac{R}{12}g+W^-
\end{array} \right)=\left( \begin{array}{cc}
a_i+b_i & 0 \\
0 & a_i-b_i
\end{array} \right).
\end{equation}
Therefore for Einstein four-manifolds, $\mathfrak{R}$ is half two-positive if and only if $(a_1+a_2)\pm(b_1+b_2)>0$.

Huisken \cite{Huiskin} observed that,

\begin{lemma} \label{WeylBergerdecomp}
Berger curvature decomposition works for every algebraic curvature tensor with constant trace on four-manifolds.
\end{lemma}
\Pf. The proof is the same as Berger's proof of Proposition \ref{Bergerdecomposition}. For the reader's convenience, we provide the proof in the appendix.

\

Let us first provide an alternative proof of the Weitzenbo\"ock formula for Einstein metrics on four-manifolds. 

\Pf of Theorem \ref{FourWeitzenbock}. The proof follows from Section 2.3 of the author's thesis \cite{Wuthesis}. Recall that in Lemma 7.2 of \cite{Hamilton}, Hamilton proved that for an Einstein manifold $(M^n, g)$ with $\mathrm{Ric}=\lambda g$,

\begin{equation} \label{HamEquation}
\Delta R_{ijkl}+2Q(\mathrm{Rm})_{ijkl}=2\lambda R_{ijkl},
\end{equation}
where $Q(\mathrm{Rm})_{ijkl}=B(\mathrm{Rm})_{ijkl}-B(\mathrm{Rm})_{ijlk}+B(\mathrm{Rm})_{ikjl}-B(\mathrm{Rm})_{iljk}$, and $B(\mathrm{Rm})_{ijkl}=g^{mn}g^{pq}R_{imjp}R_{knlq}$.

For an Einstein four-manifold, by the standard curvature decomposition, $W=\mathrm{Rm}-\frac{\lambda}{6}g\circ g$. Since $W$ is traceless, we get

\begin{equation*}
\begin{split}
\Delta|W|^2=&2\langle \Delta W, W\rangle+2|\nabla W|^2\\
=&2\langle \Delta\mathrm{Rm}, W\rangle+2|\nabla W|^2\\
=&4\langle \lambda\mathrm{Rm}-Q(\mathrm{Rm}), W\rangle+2|\nabla W|^2\\
=&4\lambda|W|^2-4\langle Q(\mathrm{Rm}), W\rangle+2|\nabla W|^2.
\end{split}
\end{equation*}

For self-dual and anti-self-dual Weyl curvature, similarly we get

\begin{equation*}
\begin{split}
\Delta|W^{\pm}|^2
=&4\lambda|W^{\pm}|^2-4\langle Q(\mathrm{Rm})^{\pm}, W^{\pm}\rangle+2|\nabla W^{\pm}|^2.
\end{split}
\end{equation*}

Using Berger curvature decomposition, it is a direct computation that (see \cite{Wuthesis})
\begin{equation*}
\begin{split}
Q(\mathrm{Rm})_{1212}=Q(\mathrm{Rm})_{3434}=&a_1^2+b_1^2+2a_2a_3+2b_2b_3,\\
Q(\mathrm{Rm})_{1234}=&2a_1b_1+2a_2b_3+2a_3b_2,\\
Q(\mathrm{Rm})_{1313}=Q(\mathrm{Rm})_{2424}=&a_2^2+b_2^2+2a_1a_3+2b_1b_3,\\
Q(\mathrm{Rm})_{1342}=&2a_2b_2+2a_1b_3+2a_3b_1,\\
Q(\mathrm{Rm})_{1414}=Q(\mathrm{Rm})_{2323}=&a_3^2+b_3^2+2a_1a_2+2b_1b_2,\\
Q(\mathrm{Rm})_{1423}=&2a_3b_3+2a_1b_2+2a_2b_1,\\
Q(\mathrm{Rm})_{ijik}=&0, \text{ if }i\neq j\neq k.
\end{split}
\end{equation*}
Recall that $B$ has symmetries $B_{ijkl}=B_{jilk}=B_{klij}$, so we compute

\begin{equation*}
\begin{split}
\langle Q(\mathrm{Rm})^{\pm},W^{\pm}\rangle=9(\bar a_1\pm b_1)(\bar a_2\pm b_2)(\bar a_3\pm b_3)=9\det W^{\pm},
\end{split}
\end{equation*}
where $\bar a_i=a_i-\frac{\lambda}{3}$, which finishes the proof.

\

Now using the Weitzenb\"ock formula and Berger curvature decomposition, we classify Einstein four-manifolds of half two-nonnegative curvature operator.

\Pf of Theorem \ref{half}. We follow the arguments in \cite{Wuthesis}.
Without loss of generality, we assume $\mathfrak{R}^+$ is half two-nonnegative.

If $R=0$, then $\mathfrak{R}^+=W^+$. Hence $\mathfrak{R}^+$ two-nonnegative implies $W^+\equiv 0$.

If $R>0$, without loss of generality, we assume $\text{Ric}=g$. Integrate the Weitzenb\"ock formula,
\begin{equation} \label{Weitzenbock}
0 = \int_M \Delta |W^+|^2 = \int_M 2|\nabla W^+|^2+(4|W^+|^2-36\det W^+).
\end{equation}

Denote $a=a_1+b_1-\frac{1}{3}, b=a_2+b_2-\frac{1}{3}, c=a_3+b_3-\frac{1}{3}$ be eigenvalues of $W^+$. By Berger curvature decomposition we have
\begin{equation*}
f = 4|W^+|^2-36\det W^+=8(a^2+ac+c^2)+36ac(a+c),
\end{equation*}

(1). If $\mathfrak{R}^+$ is two-positive, then $0\leq c<\frac{2}{3}$, $-2c\leq a\leq -\frac{c}{2}$. Taking the first derivative of $f$, we get $f_a=(2a+c)(8+36c)\leq 0$. Hence the minimum of $f$ is attained at $a=-\frac{c}{2}$, at which
\begin{equation} \label{f}
f = 6c^2(2-3c)\geq 0,
\end{equation}
with equality if and only if $c=0$, i.e. $W^+=0$.

Therefore by equation (\ref{Weitzenbock}) we get $W^+\equiv 0$, i.e. $(M,\ g)$ is anti-self-dual.

\

(2). If $\mathfrak{R}^+$ is two-nonnegative and $W^+\not \equiv 0$, then by (\ref{Weitzenbock}) and (\ref{f}) we get
\begin{equation} \label{7}
\nabla W^+=0,\hspace{1cm} a=b=-\frac{1}{3},\ c=\frac{2}{3}.
\end{equation}
Therefore by a theorem of Derdzinski \cite{Der}, $(M,\ g)$ is a K\"ahler-Einstein manifold.

\

(3). Denote $x=a_1-b_1-\frac{1}{3}, y=a_2-b_2-\frac{1}{3}, z=a_3-b_3-\frac{1}{3}$ be eigenvalues of $W^-$. If $\mathfrak{R}^+$ is two-nonnegative and $\mathfrak{R}$ is four-nonnegative, assuming $W^+\not \equiv 0$, then by equation \eqref{7},
\begin{equation*}
a+b+x+y+\frac{4}{3}=x+y+\frac{2}{3}\geq 0,
\end{equation*}
so $\mathfrak{R}^-$ is also two-nonnegative. By the same argument as above, if $W^-\not\equiv 0$, then
\begin{equation*}
\nabla W^-=0,\hspace{1cm} x=y=-\frac{1}{3},\ z=\frac{2}{3}.
\end{equation*}
Therefore $\nabla R=0$ (hence $(M,g)$ is locally symmetric) and $\mathfrak{R}$ has eigenvalues $\{0, 0, 1, 0, 0, 1\}$. By the classification of four-dimensional symmetric spaces, it is isometric to $(S^2\times S^2, g_0\oplus g_0)$ or its finite quotient.


\section{Proof of Theorem 1.3 and Theorem 1.4}

Similar to the proof of Theorem \ref{FourWeitzenbock}, we first derive a Weitzenb\"ock formula for the curvature tensor. We start from the following basic lemma (see \cite{CSW}),

\begin{lemma} \label{lemma}
Let $(M^n, g)$ be a generalized quasi-Einstein manifold with $\mathrm{Ric}_f^m=\lambda g$, then
\begin{equation*}
\begin{split}
\nabla_i R_{jk}-\nabla_j R_{ik} =& (\nabla_i\lambda g_{jk}-\nabla_j\lambda g_{ik})-R_{ijkl}\nabla_l f\\
+&\frac{1}{m}(\lambda g_{ik}\nabla_j f-\lambda g_{jk}\nabla_i f+R_{jk}\nabla_i f-R_{ik}\nabla_j f),\\
\nabla_i R =& (n-1)\nabla_i\lambda+2R_{ij}\nabla_j f-\frac{2}{m}R_{ij}\nabla_j f+\frac{2}{m}[R-(n-1)\lambda]\nabla_i f.
\end{split}
\end{equation*}
\end{lemma}

\Pf. It follows directly from the Ricci identity,
\begin{equation*}
\begin{split}
&\nabla_i R_{jk}-\nabla_j R_{ik}\\
=& \nabla_i(\lambda g_{jk}+\frac{1}{m}\nabla_j f\nabla_k f-\nabla_j\nabla_k f)
-\nabla_j(\lambda g_{ik}+\frac{1}{m}\nabla_i f\nabla_k f-\nabla_i\nabla_k f)\\
=& (\nabla_i\lambda g_{jk}-\nabla_j\lambda g_{ik})+\frac{1}{m}(\nabla_i\nabla_k f\nabla_j f-\nabla_j\nabla_k f\nabla_i f)
+(\nabla_j\nabla_i\nabla_k f-\nabla_i\nabla_j\nabla_k f)\\
=& (\nabla_i\lambda g_{jk}-\nabla_j\lambda g_{ik})-R_{ijkl}\nabla_l f+\frac{1}{m}(\nabla_i\nabla_k f\nabla_j f-\nabla_j\nabla_k f\nabla_i f)\\
=& (\nabla_i\lambda g_{jk}-\nabla_j\lambda g_{ik})-R_{ijkl}\nabla_l f+\frac{1}{m}(\lambda g_{ik}\nabla_j f-\lambda g_{jk}\nabla_i f
+R_{jk}\nabla_i f-R_{ik}\nabla_j f).
\end{split}
\end{equation*}

Taking the trace we get the second equation.

\qed

Using Lemma \ref{lemma} and Hamilton's argument we get,

\begin{prop} \label{Laplacian}
Let $(M^n, g)$ be a generalized quasi-Einstein manifold with $\mathrm{Ric}_f^m=\lambda g$, then
\begin{equation*}
\begin{split}
\Delta_f R_{ijkl} =& 2\lambda R_{ijkl}-2Q(R)_{ijkl}+(\nabla^2\lambda\circ g)_{ijkl} +\frac{1}{m}[(\text{Ric}-\lambda g)\circ\nabla^2 f]_{ijkl}\\
+&\frac{1}{m^2}[(\text{Ric}-\lambda g)\circ df\otimes df]_{ijkl}+\frac{1}{m}(\nabla\lambda\otimes\nabla f\circ g)_{ijkl}\\
+&\frac{1}{m}[R_{ipkl}\nabla_j f\nabla^p f+R_{ijkp}\nabla_lf\nabla^p f-R_{jpkl}\nabla_i f\nabla^p f-R_{ijlp}\nabla_k f\nabla^p f].
\end{split}
\end{equation*}
\end{prop}

\Pf. By the Ricci identity, we get (see Lemma 7.2 in Hamilton \cite{Hamilton})
\begin{equation*}
\begin{split}
\Delta R_{ijkl}=&\nabla^p\nabla_p R_{ijkl}\\
=&\nabla^p\nabla_iR_{pjkl}-\nabla_p\nabla_j R_{pikl}\\
=& \nabla_i(\nabla_k R_{jl}-\nabla_l R_{jk})-\nabla_j(\nabla_k R_{il}-\nabla_l R_{ik})\\
-&2Q(R)_{ijkl}+(R_i^{\ q}R_{qjkl}-R_j^{\ q}R_{qikl})
\end{split}
\end{equation*}

Applying Lemma \ref{lemma} repeatedly to the first two terms on the right hand side, we have
\begin{equation*}
\begin{split}
& \nabla_i(\nabla_k R_{jl}-\nabla_l R_{jk})-\nabla_j(\nabla_k R_{il}-\nabla_l R_{ik})\\
=&\nabla_i[(\nabla_k\lambda g_{jl}-\nabla_l\lambda g_{jk})-R_{kljp}\nabla^p f\\
&+\frac{1}{m}(\lambda g_{jk}\nabla_lf-\lambda g_{jl}\nabla_k f+R_{jl}\nabla_k f-R_{jk}\nabla_l f)]\\
-&\nabla_j[(\nabla_k\lambda g_{il}-\nabla_l\lambda g_{ik})-R_{klip}\nabla^p f\\
&+\frac{1}{m}(\lambda g_{ik}\nabla_lf-\lambda g_{il}\nabla_k f+R_{il}\nabla_k f-R_{ik}\nabla_l f)]\\
=& (\nabla^2\lambda\circ g)_{ijkl}-(\nabla_i R_{jpkl}\nabla^p f +\nabla_j R_{pikl}\nabla^p f)
+(R_{ipkl}\nabla_j\nabla^p f-R_{jpkl}\nabla_i\nabla^p f)\\
&+\frac{1}{m}[(\text{Ric}-\lambda g)\circ\nabla^2 f]_{ijkl}+\frac{1}{m}[(\nabla_iR_{jl}-\nabla_j R_{il})\nabla_k
f-(\nabla_i R_{jk}-\nabla_j R_{ik})\nabla_l f]\\
=& (\nabla^2\lambda\circ g)_{ijkl}+\nabla_p R_{ijkl}\nabla^p f+(R_{ipkl}\nabla_j\nabla^p f-R_{jpkl}\nabla_i\nabla^p f)\\
&+\frac{1}{m}(R_{kpij}\nabla_l f\nabla^p f-R_{lpij}\nabla_k f\nabla^p f) +\frac{1}{m}[(\text{Ric}-\lambda g)\circ\nabla^2 f]_{ijkl}\\
&+\frac{1}{m}(\nabla\lambda\otimes\nabla f\circ g)_{ijkl} +\frac{1}{m^2}[(\text{Ric}-\lambda g)\circ df\otimes df]_{ijkl}.
\end{split}
\end{equation*}

Therefore we get
\begin{equation*}
\begin{split}
&\Delta_f R_{ijkl}\\
=& -2Q(\mathrm{Rm})_{ijkl}+(R_{ipkl}R_j^{\ p}-R_{jpkl}R_i^{\ p})+(\nabla^2\lambda\circ g)_{ijkl}+(R_{ipkl}\nabla_j\nabla^p f-R_{jpkl}\nabla_i\nabla^p f) \\
&+\frac{1}{m}(R_{kpij}\nabla_l f\nabla^p f-R_{lpij}\nabla_k f\nabla^p f) +\frac{1}{m}[(\text{Ric}-\lambda g)\circ\nabla^2 f]_{ijkl}\\
&+\frac{1}{m}(\nabla\lambda\otimes\nabla f\circ g)_{ijkl}+\frac{1}{m^2}[(\text{Ric}-\lambda g)\circ df\otimes df]_{ijkl}\\
=&2\lambda R_{ijkl}-2Q(\mathrm{Rm})_{ijkl}+(\nabla^2\lambda\circ g)_{ijkl}+\frac{1}{m}(\nabla\lambda\otimes\nabla f\circ g)_{ijkl}\\
&+\frac{1}{m}[(\text{Ric}-\lambda g)\circ\nabla^2 f]_{ijkl} +\frac{1}{m^2}[(\text{Ric}-\lambda g)\circ df\otimes df]_{ijkl}\\
&+\frac{1}{m}[R_{ipkl}\nabla_j f\nabla^p f-R_{jpkl}\nabla_i f\nabla^p f+R_{kpij}\nabla_l f\nabla^p f-R_{lpij}\nabla_k f\nabla^p f].
\end{split}
\end{equation*}

\qed

Applying the standard curvature decomposition and Berger curvature decomposition, we prove the Weitzenb\"ock for $W^{\pm}$.

\Pf of Theorem \ref{UnifiedWeitzenbock}. We need to express $\Delta_f R_{ijkl}$ in terms of Weyl curvature using the standard curvature decomposition,
\begin{equation*}
\begin{split}
R_{ijkl} &= -\frac{R}{2(n-1)(n-2)}(g\circ g)_{ijkl}+\frac{1}{n-2}(Ric\circ g)_{ijkl}+W_{ijkl}.
\end{split}
\end{equation*}
First we have (see Catino and Mantegazza \cite{CaMa})
\begin{equation*}
\begin{split}
2Q(\mathrm{Rm})_{ijkl} =&2Q(W)_{ijkl}+\frac{2(n-1)|\text{Ric}|^2-2R^2}{2(n-1)(n-2)^2}(g\circ g)_{ijkl} +\frac{2}{n-2}(R_{ik}R_{lj}-R_{il}R_{jk})\\
-&\frac{2}{(n-2)^2}(\text{Ric}^2\circ g)_{ijkl}+\frac{2R}{(n-1)(n-2)^2}(\text{Ric}\circ g)_{ijkl}\\
+&\frac{2}{n-2}(W_{ipkq}R^{pq}g_{jl}-W_{jpkq}R^{pq}g_{il}+W_{jplq}R^{pq}g_{ik}-W_{iplq}R^{pq}g_{jk}).
\end{split}
\end{equation*}
Similarly we compute
\begin{equation*}
\begin{split}
&\frac{1}{m}[R_{ipkl}\nabla_j f\nabla^p f-R_{jpkl}\nabla_i f\nabla^p f+R_{kpij}\nabla_l f\nabla^p f-R_{lpij}\nabla_k f\nabla^p f]\\
=& \frac{1}{m}[W_{ipkl}\nabla_j f\nabla^p f+W_{ijkp}\nabla_l f\nabla^p f-W_{jpkl}\nabla_i f\nabla^p f-W_{ijlp}\nabla_k f\nabla^p f]\\
-&\frac{2R}{m(n-1)(n-2)}(g\circ df\otimes df)_{ijkl}\\
+&\frac{2}{m(n-2)}\Big[(\text{Ric}\circ df\otimes df)_{ijkl}+[\nabla f\otimes\text{Ric}(\nabla f)\circ g]_{ijkl}\Big].
\end{split}
\end{equation*}

\

Since $W$ is traceless, $\langle\alpha\circ g,W\rangle=0$ for any $(0,2)$-tensor $\alpha$. Therefore we get the following Weitzenb\"ock formula for Weyl curvature,
\begin{equation} \label{formula}
\begin{split}
&\Delta_f |W|^2\\
=& 2\langle\Delta_f W, W\rangle+2|\nabla W|^2\\
=& 2\langle\Delta_f\mathrm{Rm}, W\rangle+2|\nabla W|^2\\
=&2|\nabla W|^2+4\lambda |W|^2-W^{ijkl}Q(W)_{ijkl}-\frac{1}{2(n-2)}(\text{Ric}\circ\text{Ric})_{ijkl}W^{ijkl}\\
+&\frac{1}{2m}W^{ijkl}[W_{ipkl}\nabla_j f\nabla^p f+W_{ijkp}\nabla_l f\nabla^p f-W_{jpkl}\nabla_i f\nabla^p f-W_{ijlp}\nabla_k f\nabla^p f]\\
+&\frac{1}{2m}(\text{Ric}\circ \nabla^2 f)_{ijkl}W^{ijkl}+\Big[\frac{1}{2m^2}+\frac{1}{m(n-2)}\Big](\text{Ric}\circ df\otimes df)_{ijkl}W^{ijkl}\\
=&2|\nabla W|^2+4\lambda |W|^2-4\langle W, Q(W)\rangle-\frac{2}{n-2}\langle\text{Ric}\circ\text{Ric},W\rangle\\
+&\frac{8}{m}|\iota_{\nabla f}W|^2+\frac{2}{m}\langle\text{Ric}\circ \nabla^2 f,W\rangle+\Big[\frac{2}{m^2}+\frac{4}{m(n-2)}\Big]\langle\text{Ric}\circ
df\otimes df,W\rangle.
\end{split}
\end{equation}

\

Next we derive the Weitzenb\"ock formula for $W^{\pm}$ on four-manifolds. For any $\sigma\in S_4$, denote $\sigma(ijkl)=(\sigma_i\sigma_j\sigma_k\sigma_l)$, it is a direct computation that
\begin{equation*}
\begin{split}
\langle(\alpha\circ g)_{ijkl}, W^{\sigma_i\sigma_j\sigma_k\sigma_l}\rangle
=&\frac{1}{4}(\alpha_{ik}g_{jl}+\alpha_{jl}g_{ik}-\alpha_{il}g_{jk}-\alpha_{jk}g_{il})W^{\sigma_i\sigma_j\sigma_k\sigma_l}=0,
\end{split}
\end{equation*}
that is
\begin{equation*}
\begin{split}
\langle(\alpha\circ g)^{\pm}, W^{\pm}\rangle=0.
\end{split}
\end{equation*}

Recall that for any $(ij)$, $(i'j')$ is defined to be the pair such that $e_i\wedge e_j\pm e_{i'}\wedge e_{j'}\in \wedge^{\pm} M$, and for any $(0,4)$-tensor T,
\begin{equation*}
\begin{split}
T^{\pm}_{ijkl}=&\frac{1}{4} T(e_i\wedge e_j\pm e_{i'}\wedge e_{j'},e_k\wedge e_l\pm e_{k'}\wedge e_{l'})\\
=&\frac{1}{4}(T_{ijkl}\pm T_{ijk'l'}\pm T_{i'j'kl}+T_{i'j'k'l'}).
\end{split}
\end{equation*}

So we get
\begin{equation*}
\begin{split}
&\Delta|W^{\pm}|^2\\
=&2\langle\Delta W^{\pm},W^{\pm}\rangle+2|\nabla W^{\pm}|^2\\
=&2\langle\Delta \mathrm{Rm}^{\pm},W^{\pm}\rangle+2|\nabla W^{\pm}|^2\\
=&\frac{1}{8}\langle\Delta(R_{ijkl}\pm R_{ijk'l'}\pm R_{i'j'kl}+R_{i'j'k'l'}),(W^{ijkl}\pm W^{ijk'l'}\pm W^{i'j'kl}+W^{i'j'k'l'})\rangle\\
&+2|\nabla W^{\pm}|^2.
\end{split}
\end{equation*}

Same as in the proof of Theorem \ref{FourWeitzenbock}, using the Berger curvature decomposition for $W^{\pm}$, we get
\begin{equation*}
W^{\pm\,ijkl}Q(W)_{ijkl}^{\pm}=36\det W^{\pm}.
\end{equation*}
Therefore from equation \eqref{formula}, we have
\begin{equation*}
\begin{split}
&\Delta_f |W^{\pm}|^2\\
=& 2|\nabla W^{\pm}|^2+4\lambda |W^{\pm}|^2-36\det W^{\pm} -\langle(\text{Ric}\circ\text{Ric})^{\pm},W^{\pm}\rangle\\
+&\frac{2}{m}\langle(\text{Ric}\circ \nabla^2 f)^{\pm},W^{\pm}\rangle +(\frac{2}{m}+\frac{2}{m^2})\langle(\text{Ric}\circ df\otimes
df)^{\pm},W^{\pm}\rangle\\
+&\frac{1}{2m}W^{\pm\,ijkl}\Big[W_{ipkl}\nabla_j f\nabla^p f+W_{ijkp}\nabla_l f\nabla^p f-W_{jpkl}\nabla_i f\nabla^p f
-W_{ijlp}\nabla_k f\nabla^p f\Big]^{\pm}\\
=& 2|\nabla W^{\pm}|^2+4\lambda |W^{\pm}|^2-36\det W^{\pm}\\
+&(1+\frac{2}{m})\langle(\text{Ric}\circ\text{Ric})^{\pm},W^{\pm}\rangle +(2+\frac{4}{m})\langle(\text{Ric}\circ \nabla^2 f)^{\pm},W^{\pm}\rangle\\
+&\frac{1}{2m}W^{\pm\,ijkl}\Big[W_{ipkl}\nabla_j f\nabla^p f\pm W_{i'pkl}\nabla_{j'} f\nabla^p f\pm W_{ipk'l'}\nabla_jf\nabla^pf
+W_{i'pk'l'}\nabla_{j'}f \nabla^pf\Big].
\end{split}
\end{equation*}

By the symmetry $W^{\pm}_{ijkl}=\pm W^{\pm}_{i'j'kl}=\pm W^{\pm}_{ijk'l'}=W^{\pm}_{i'j'k'l'}$, we have

\begin{equation*}
\begin{split}
&\frac{1}{2m}W^{\pm\,ijkl}[W_{ipkl}\nabla_j f\nabla^p f\pm W_{i'pkl}\nabla_{j'} f\nabla^p f\pm W_{ipk'l'}\nabla_j f\nabla^p f
+W_{i'pk'l'}\nabla_{j'} f\nabla^p f]\\
=&\frac{2}{m}W^{\pm\,ijkl}W_{ipkl}\nabla_j f\nabla^p f\\
=&\frac{8}{m}\langle\iota_{\nabla f}W^{\pm},\iota_{\nabla f}W\rangle.
\end{split}
\end{equation*}

By definition $\langle W^{+}, W^-\rangle=0$. Using Berger curvature decomposition, it is easy to verify that
\begin{lemma}
Let $(M,\ g)$ be a four-manifold. Then for any $f\in C^{\infty}(M)$,
\begin{equation*}
\begin{split}
\langle\iota_{\nabla f}W^+,\iota_{\nabla f}W^-\rangle=\frac{1}{4}W^{+\,ipkl}W^-_{iqkl}\nabla_p f\nabla^q f &=0,\\
|\iota_{\nabla f}W^{\pm}|^2=\frac{1}{4}W^{\pm\,ipkl}W^{\pm}_{iqkl}\nabla_p f\nabla^q f
&=\frac{1}{4}|W^{\pm}|^2|\nabla f|^2.
\end{split}
\end{equation*}
\end{lemma}

Therefore we obtain
\begin{equation*}
\begin{split}
\frac{8}{m}\langle\iota_{\nabla f}W^{\pm},\iota_{\nabla f}W\rangle =\frac{8}{m}\langle\iota_{\nabla f}W^{\pm},\iota_{\nabla f}(W^+ +W^-)\rangle
=\frac{2}{m}|W^{\pm}|^2|\nabla f|^2.
\end{split}
\end{equation*}
This completes the proof of Theorem \ref{UnifiedWeitzenbock}.

\qed

Similarly to Einstein four-manifolds, we observe that for any four-manifold of half two-nonnegative curvature operator,
\begin{lemma} \label{halfPIC}
Let $(M^4,\ g)$ be a four-manifold. If $\mathfrak{R}^{\pm}$ is two-nonnegative, then
\begin{equation*}
R|W^{\pm}|^2-36\det W^{\pm}\geq 0.
\end{equation*}

Furthermore, if $\mathfrak{R}^{\pm}$ is two-positive, the equality holds if and only if $W^{\pm}=0$. If $\mathfrak{R}^{\pm}$ is two-nonnegative, the
equality holds if and only if $W^{\pm}=0$, or $W^{\pm}$ has eigenvalues $\{-\frac{R}{12},\ -\frac{R}{12},\ \frac{R}{6}\}$.
\end{lemma}

\Pf of Lemma \ref{halfPIC}. The proof is similar to Theorem \ref{half}. Without loss of generality, we assume $\mathfrak{R}^-$ is two-nonnegative. Denote $x\leq y\leq z$ be eigenvalues of $\mathfrak{R}^-$. Let
\begin{equation*}
f=R|W^-|^2-36\det W^-=2R(x^2+xz+z^2)+36xz(x+z).
\end{equation*}
Taking the derivative, we get $f_x =2(2x+z)(R+18z)\leq 0$, so the minimum of $f$ is attained at $x=-\frac{z}{2}$, at which
\begin{equation*}
f =\frac{3}{2}z^2(R-6z).
\end{equation*}
If $\mathfrak{R}^-$ is two-nonnegative, then $x+y+\frac{R}{6}\geq0$, so $\frac{R}{6}-z\geq0$. Therefore $f\geq 0$.

Moreover, if $\mathfrak{R}^-$ is two-positive, $f=0$ if and only if $W^-=0$. If $\mathfrak{R}^-$ is two-nonnegative, $f=0$ if and only if $W^-=0$, or $W^-$ has eigenvalues $\{-\frac{R}{12},\ -\frac{R}{12},\ \frac{R}{6}\}$.

\qed

\Pf of Theorem \ref{harmonicWeyl}. Assume for example $W^+$ is harmonic and $\mathfrak{R}^+$ is two-nonnegative. By the integral Weitzenb\"ock formula \eqref{intWeitzenbock}, we get
\begin{equation*}
\begin{split}
0=& \int_M 2|\nabla W^+|^2+R|W^+|^2-36\det W^+.
\end{split}
\end{equation*}
Therefore by Lemma \ref{halfPIC}, $$\nabla W^+\equiv 0,\ \ \mathrm{and}\ \ R|W^+|^2-36\det W^+\equiv 0.$$

If $\mathfrak{R}^+$ is two-positive, then $W^+\equiv 0$.

If $\mathfrak{R}^+$ is two-nonnegative, then by Lemma \ref{halfPIC}, either $W^+\equiv 0$ or $W^+$ has eigenvalues $\{-\frac{R}{12},\ -\frac{R}{12},\
\frac{R}{6}\}$. If $W^+\not\equiv 0$ then $\nabla W^+\equiv 0$ implies that $R\equiv$const, then by a theorem of Derdzinski (Proposition 5 in \cite{Der}), $g$ is a K\"ahler metric.
\vspace{0.4cm}

If in addition $(M,g,f)$ is a gradient shrinking Ricci soliton.  If $W^{\pm}=0$, then by the work of Chen-Wang \cite{ChenWang} or Cao-Chen \cite{CaoChen}, $(M,g,f)$ must be isometric to $(S^4,\ g_0)$ or $(\mathbb{C}P^2,\ g_{FS})$.

If $W^{\pm}\not\equiv 0$, then by the soliton equation $R+\Delta f=4\lambda$, $R\equiv$const implies $f\equiv$const, therefore $g$ is a K\"ahler-Einstein metric.
\qed

\


\section{Appendix: Proof of Berger curvature decomposition}

The proof is translated directly from Berger's paper \cite{Berger}.

\

Let $P\subset T_p M$ be the 2-plane such that the sectional curvature attains its minimum on $P$. Let $P^{\bot}$ be the two plane that is orthogonal to $P$, choose $e_1\in P,\ e_2\in P^{\bot}$ such that $K(e_1, e_2)\geq K(X,Y)$ for any $X\in P,\ Y\in P^{\bot}$. Expand $\{e_1,\ e_2\}$ to an orthonormal basis $\{e_1,\ e_2,\ e_3,\ e_4\}$ such that $e_3\in P,\ e_4\in P^{\bot}$.

By the choice of $P$, $K(X,Y)\geq K(e_1,e_3)$ for any $X,\ Y\in T_p M$, in particular, $K(X,e_3)\geq K(e_1,e_3)$, $K(Y,e_1)\geq K(e_1,e_3)$ for any $X,\ Y\in T_p M$.

Let $X=e_1\cos t +e_2\sin t$, $Y=e_3\cos t +e_2\sin t$, $0\leq t\leq\delta$, by variation principle we get
\begin{equation*}
\begin{split}
0 = \frac{d}{dt}\Big|_{t=0} K(X,e_3) &= 2R_{1323},\ \ \ \ 0 = \frac{d}{dt}\Big|_{t=0} K(Y,e_1) = 2R_{1312},
\end{split}
\end{equation*}

Similarly let $X=e_1\cos t +e_4\sin t$, $Y=e_3\cos t +e_4\sin t$, we get $R_{1343}=0$, $R_{1314}=0$.

By Lemma \ref{BergerLemma}, $K(e_2, e_4)=K(e_1,e_3)$, so by the same argument as above, we have $R_{2124}=R_{2324}=R_{4142}=R_{4243}=0$.
\vspace{0.4cm}

On the other hand, we have $K(e_1, e_2)\geq K(X,Y)$ for any $X\in P,\ Y\in P^{\bot}$, in particular, $K(e_1, e_2)\geq K(e_1,X)$ for any $X\in P^{\bot}$. Let $X=e_2\cos t +e_4\sin t$ by variation principle we have $R_{1214}=0$; also $K(e_1, e_2)\geq K(X, e_2)$ for any $X\in P$. Let $X=e_1\cos t +e_3\sin t$ we get $R_{2123}=0$.

Again by Lemma \ref{BergerLemma}, $K(e_3, e_4)=K(e_1, e_2)\geq K(X,Y)$ for any $X\in P,\ Y\in P^{\bot}$, so we get $R_{3432}=R_{4341}=0$. Therefore we proved (1) and (2).
\vspace{0.4cm}

Since $K(X,Y)\geq K(e_1,e_3)$ for any $X,\ Y\in T_p M$, we obtain $K(ae_1+be_2,ce_3+de_4)\geq K(e_1,e_3)$ for any $a,b,c,d$ such that $a^2+b^2=1,\ c^2+d^2=1$.

Choose $a=b=c=d=\frac{1}{\sqrt{2}}$, we get
\begin{equation*}
K(e_1,e_3) \leq K(\frac{1}{\sqrt{2}}e_1 +\frac{1}{\sqrt{2}}e_2,
\frac{1}{\sqrt{2}}e_3 +\frac{1}{\sqrt{2}}e_4) = \frac{1}{2}[R_{1313}+R_{1414} +R_{1324}+R_{1423}]
\end{equation*}
Therefore $$R_{1423}-R_{1342}\geq R_{1313}-R_{1414}.$$

Similarly choosing $a=b=c=\frac{1}{\sqrt{2}},\ d=-\frac{1}{\sqrt{2}}$ we get $R_{1342}-R_{1423}\geq R_{1313}-R_{1414}$, therefore $$|R_{1342}-R_{1423}|\leq
R_{1414}-R_{1313}.$$

Apply the same argument to $K(ae_1+be_3,ce_2+de_4)\leq K(e_1,e_2)$, and $K(ae_1+be_4,ce_2+de_3)\geq K(e_1,e_3)$, we get
\begin{equation*}
\begin{split}
|R_{1234}-R_{1342}|\leq& R_{1313}-R_{1212},\\
|R_{1423}-R_{1234}|\leq& R_{1414}-R_{1212}.
\end{split}
\end{equation*}

\qed

\end{document}